\documentclass[titlepagea4paper,11pt]{article}
\usepackage[arrow,matrix]{xy}
\usepackage{amssymb}
\usepackage{amsbsy,amsmath}
\newtheorem{thm}{Theorem}[section]
\newtheorem{lemma}[thm]{Lemma}

\newtheorem{rem}[thm]{Remark}
\newcommand{\bb}{\mathbb}
\newcommand{\ratmap}{\cdots\!\rightarrow}
\newcommand{\seq}{\longrightarrow}

\newcommand{\can}{-K_{{\bb P}({\cal E})}}
\newenvironment{proof}{\vspace{3pt}\indent
                       \textsc{\bf Proof:}\quad }
                       {\hfill$\square$\vspace{3pt}}
\newtheorem{example}[thm]{Example}
\newtheorem{construction}[thm]{Construction}

\newenvironment{constr}{\begin{construction}\rm}{\hfill$\bigtriangleup$\end{construction}}

\newenvironment{remark}{\begin{rem}\rm}{\end{rem}}

\title{A note on octic hypersurfaces with many nodes}
\author{Marco K\"uhnel}
\date{\today}

\begin{document}
\maketitle
\abstract{In the present note we use rank-2-bundles over ${\bb P}^3$ to construct octic hypersurfaces with many nodes. We give an example with 128 nodes.}

\section*{Introduction}

Given a nodal octic $S$ in ${\bb P}^3$, one can construct a Calabi-Yau-threefold, which is a desingularization of a double cover of ${\bb P}^3$ ramified over
$S$ (cf. \cite{cl83}). So the usual way is first to prove existence of octics with a certain number of nodes and then to construct 
from this Calabi-Yau-threefolds with certain Euler numbers. 
This can also be done for octic arrangements with isolated singularities (cf. \cite{cy99,cs99}). 

It has been proved, that for every given number not larger than $108$ there are octic hypersurfaces with this number of nodes 
(\cite{we87,bo90}). Furthermore
there are examples for many numbers between 108 and 168. For the probably most complete list of constructed numbers we refer to \cite{labs}.
In the references the interested reader may find the most important sources for constructions of octics with many nodes. 
The author likes to thank Oliver Labs for pointing out several references.
Miyaoka proved in \cite{mi84} an upper bound of $174$ possible nodes. 

In this note the way is reversed: Calabi-Yau-threefolds, which are desingularizations of double covers of ${\bb P}^3$ are constructed, what gives rise to
octics with a certain number of nodes. To construct the Calabi-Yau-threefolds we look at hypersurfaces of projectivized rank-2-bundles over ${\bb P}^3$.
For this class it turns out that, indeed, the octics can only have nodes as singularities. Finally, 
we construct an example with $128$ nodes as our maximal case up to now.

\section{Construction of the octic hypersurface}

Let ${\cal E}\stackrel{p}{\seq}{\bf P}^3$ be a rank-2-bundle and  $s\in H^0(\can)$ a section such that $X=\{s=0\}$ is smooth. 
We denote $\gamma({\cal E}):=\deg(c_1^2({\cal E})-4c_2({\cal E})).$ This is invariant under tensorizing ${\cal E}$ with line bundles,
and moreover
$$c_3(X)=-8\gamma-168,$$
what can be computed by standard methods. Hence $\gamma$ is a topological invariant of $X$.

The restriction to $X$ of the projection $p$ 
we call $\pi:X\seq{\bb P}^3$. The map $\pi$ is a generic double cover. If we take the Stein factorization
$$X\stackrel{c}{\seq}X'\stackrel{\phi}{\seq}{\bb P^3},\quad\quad \pi=\phi\circ c,$$
the map $\phi$ is a double cover, whose ramification  locus we denote by $B\subset{\bb P}^3$. 
We call $B$ also the ramification divisor of $\pi$ and denote $\tilde B:=\pi^{-1}(B)$.
Another description of $B$ is obtained by looking at the discriminant map.

\begin{constr}\label{3discr}
Let $X=\{s=0\}$, with $s\in H^0(-K_{{\bf P}({\cal E})})$ and
$$B:=\{p\in{\bf P}^3 | \pi \mbox{ is locally in } p\mbox{ not an \'etale covering} \}.$$
We define the discriminant
$$\Delta_{\cal E}:S^2{\cal E}\otimes{\cal F}\seq (\det({\cal E})\otimes{\cal F})^{\otimes 2}$$
by
$$\Delta_{\cal E}(\sum_{1\le i<j\le 2}c_{ij}s_is_j\otimes f):=(c_{12}^2-4c_{11}c_{22})(s_1\wedge s_2\otimes f)^{\otimes 2},$$
where $s_1,s_2$ is a ${\cal O}(U)$-basis of ${\cal E}(U)$, ${\cal F}\seq{\bf P}^3$ a line bundle and 
$f\in {\cal F}(U)$ a generator of ${\cal F}(U)$ for a small open set $U\subset{\bf P}^3$.

It is an easy computation that this definition is independent of the chosen bases.

Now we specify 
$${\cal F}=\det{\cal E}^\vee\otimes{\cal O}(4).$$
Then the discriminant is a map
$$\Delta_{\cal E}:p_*(-K_{{\bf P}({\cal E})})\seq {\cal O}(8),$$
with
$$\{\Delta_{\cal E}(p_*s)=0\}=B$$
set theoretically: in local coordinates
$$s=\sum s_{ij}x_ix_j,$$
where $[x_0:x_1]$ denotes the coordinates of the fibre, $B$ is the locus, where the zeroes of
\begin{eqnarray*}&\sum s_{ij}(z)x_ix_j=0&\label{qu}\end{eqnarray*}
are not two distinct points. By definition this is the
discriminant locus of the qudratic equation in $x_0,x_1$, given by
$$s_{01}^2-4s_{00}s_{11}=0.$$
This coincides with the discriminant locus of $p_*s$.

Since on a trivializing neighbourhood $U\subset{\bb P}^3$ the map is given by
$$\Delta_{\cal E}(t)|U=t_{12}^2-4t_{11}t_{22},$$
if $t\in H^0(p_*(-K_{{\bb P}({\cal E})}))$ and $t|U=(t_{11},t_{12},t_{22})$,
we see, that, in particular, $H^0(\Delta_{\cal E})$ is a holomorphic map.

Moreover, 
$$H^0(\Delta_{\cal E})(rt)=r^2H^0(\Delta_{\cal E})(t)$$
for $r\in{\bf C}, t\in H^0(-K_{{\bf P}({\cal E})})$. 
Hence we can projectivize, not exclude, however, that
$H^0(\Delta_{\cal E})(s')=0$ for some $s'\not= 0$. Therefore we get a rational map
$$\delta_{\cal E}:{\bf P}(H^0(-K_{{\bf P}({\cal E})}))\ratmap{\bf P}(H^0({\cal O}(8)))\cong{\bf P}^{164}.$$

Let for the moment $B':=\{z\in{\bf P}^3| H^0(\Delta_{\cal E})(s)(z)=0\}$ in the sense of ideals.
If we denote
$$P:=\{z\in{\bf P}^3 | \dim \pi^{-1}(z)=1\},$$
then we see, that
$$P=\{z\in{\bf P}^3| \sum s_{ij}x_ix_j=0 \mbox{ f\"ur alle}[x_o:x_1]\}$$
and hence
$$P=\{z\in{\bf P}^3| s_{00}(z)=s_{01}(z)=s_{11}(z)=0\}\subset Sing(B').$$
Moreover, this shows
$$P=\{z\in{\bf P}^3 | \pi^{-1}(z)\cong{\bf P}^1\}.$$

Now let $z\in Sing(B')$. If $s_{00}(z)=s_{01}(z)=s_{11}(z)=0$, then $z\in P$. 
So let us assume $s_{00}(z)\not= 0$ or $s_{01}(z)\not= 0$. Let us define
$$x:=[s_{01}(z):-2s_{00}(z)]\in p^{-1}(z).$$
Since $z\in B$, we get that $\Delta_{\cal E}(s)(z)=s_{01}(z)^2-4s_{00}(z)s_{11}(z)=0.$ Therefore
$$s(x)=s_{00}(z)s_{01}(z)^2-2s_{00}(z)s_{01}(z)^2+4s_{11}(z)s_{00}(z)^2=-s_{00}(z)\Delta_{\cal E}(s)(z)=0,$$
hence $x\in X$. 

We want to show that $x\in X$ is singular. For this we have to compute in the point $x$
\begin{eqnarray}
\frac{\partial s}{\partial x_0}=&2s_{00}x_0+s_{01}x_1&=0\label{3x1}\\
\frac{\partial s}{\partial x_1}=&s_{01}x_0+2s_{11}x_1&=0\label{3x2}\\
\frac{\partial s}{\partial z_i}=&\frac{\partial s_{00}}{\partial z_i}x_0^2+\frac{\partial s_{01}}{\partial z_i}x_0x_1+
\frac{\partial s_{11}}{\partial z_i}x_1^2&=0\label{3x3}
\end{eqnarray}
and we know moreover, since $z\in Sing(B')$, that in the point $z$ holds
\begin{eqnarray}
s_{01}^2-4s_{00}s_{11}&=&0\label{3z1}\\
2s_{01}\frac{\partial s_{01}}{\partial z_i}-4s_{11}\frac{\partial s_{00}}{\partial z_i}-4s_{00}\frac{\partial s_{11}}{\partial z_i}&=&0.
\label{3z2}
\end{eqnarray}

Using the expression for  $x$ in (\ref{3x1}), (\ref{3x2}) and (\ref{3x3}) we compute
\begin{eqnarray}
\frac{\partial s}{\partial x_0}=&2s_{00}s_{01}-2s_{00}s_{01}&=0\nonumber\\
\frac{\partial s}{\partial x_1}=&s_{01}^2-4s_{00}s_{11}&=0\nonumber\\
\frac{\partial s}{\partial z_i}=&\frac{\partial s_{00}}{\partial z_i}s_{01}^2-2\frac{\partial s_{01}}{\partial z_i}s_{00}s_{01}+
4\frac{\partial s_{11}}{\partial z_i}s_{00}^2&=\nonumber\\
=&4\frac{\partial s_{00}}{\partial z_i}s_{00}s_{11}-2\frac{\partial s_{01}}{\partial z_i}s_{00}s_{01}+
4\frac{\partial s_{11}}{\partial z_i}s_{00}^2&=\nonumber\\
=&-s_{00}(2s_{01}\frac{\partial s_{01}}{\partial z_i}-4s_{11}\frac{\partial s_{00}}{\partial z_i}-4s_{00}\frac{\partial s_{11}}{\partial z_i})&=0,\nonumber
\end{eqnarray}
with the last equation using (\ref{3z1}) as well as (\ref{3z2}).

Thus we have proved that $x\in X$ is singular. But we assumed $X$ to be smooth. Hence it is proven that
$P=Sing(B').$
In particular, $B'$ is reduced and therefore $B'=B$ in the sense of ideals.
\end{constr}

Now we know, if $X=\{s=0\}$ for some $s\in H^0(-K_{{\bf P}({\cal E})})$ then

\begin{lemma}\label{1}
$B=\delta_{\cal E}(X)\in|{\cal O}(8)|$ and $P=\{p_*s=0\}=Sing(B).$
\end{lemma}

Note that $p_*s$ gives the three local equations of $P$. For example, if ${\cal E}$ splits, then we can conclude that
$P$ is the complete intersection of three hypersurfaces of degrees $4-\sqrt\gamma,4$ and $4+\sqrt\gamma$. (Indeed, if ${\cal E}$
splits, then $\gamma$ is a square.)

Let us now specify the type of the singularities.
 
\begin{lemma}$B$ has only double points of type $A_1$ as singularities.\label{2}\end{lemma}

\begin{proof}
First, singularities of $\tilde B$ can only occur over singularities of $B$,
hence by Lemma \ref{1}
$$s_{00}(z)=s_{01}(z)=s_{11}(z)=0,$$
if $y\in\tilde B$ is singular and $z=\pi(y)$. From this we conclude again by the local descriptions that $y$ is singular in $X$. Hence $\tilde B$ is
non-singular. In particular, $B$ has only isolated singularities. 
Now we look at the rational curves $F=F_p:=\pi^{-1}(p)$ for $p\in P$. The adjunction formula yields
$${\cal O}(-2)=K_F=K_{\tilde B}|F\otimes N_{F|\tilde B}.$$
Since again by adjunction formula $\deg K_{\tilde B}|F=\tilde B.F=p^{-1}(B).F=0$ we conclude
$$N_{F|\tilde B}={\cal O}(-2)$$
and hence $p$ is a double point of type $A_1$.
\end{proof}

\begin{lemma}$|P|=64-4\gamma$.\label{3}\end{lemma}

\begin{proof}
Since by Lemma \ref{1} and Lemma \ref{2} we know that $\dim P\le 0$ and $P\cap U=\{z\in{\bb P}^3|s_{00}(z)=s_{01}(z)=s_{11}(z)=0\}$ in a trivializing neighbourhood $U$, 
we conclude
$$[P]=c_3(p_*(\can)).$$
Again by standard methods (cf. \cite[p. 423]{ha}) we compute $c_3(p_*(\can))=64-4\gamma.$
\end{proof}

This method has only limitated applications, since
by Lemma \ref{3} we see, that the number of nodes must be divisible by $4$. But there is an additional restriction:

\begin{lemma}\label{4}$\gamma({\cal E})\mod 8\in\{0,1,4\}$.\end{lemma}

\begin{proof}By definition $\gamma({\cal E})\mod 4$ is a square. The case $\gamma({\cal E})\mod 8=5$ can be excluded by the Schwarzenberger condition
$c_1({\cal E}).c_2({\cal E})\equiv 0(2)$.
\end{proof}

Since Miyaoka proved an upper bound of $174$ nodes, by Lemma \ref{4} the theoretical maximum of our method lies at $\gamma=-24$, i.e. at most $160$ nodes.
Furthermore, Lemma \ref{4} implies that the only new constructed numbers of nodes can be $124$ and $156$.

\section{Construction of some bundles}

By looking at the splitting bundles allowing for smooth $X$ we obtain immediately examples of octics with $28, 48, 60$ and $64$ nodes.
Moreover, it is not hard to see the existence of elliptic curves in ${\bf P}^3$ of degrees $d\le 7$ which are cut out by quartics
(cf. \cite{bo,hu}).
The Serre construction, more detailed explained below, then yields examples of octics with $80, 96$ and $112$ nodes.  
It is a little bit more work to see the existence of elliptic curves of degree $8$ cut out by quartics.

We use the Serre construction of rank-2-bundles to get the desired example (cf. \cite{oss}). Let $Y\subset{\bb P}^3$ be an elliptic curve of degree $8$. Then there is 
a rank-2-bundle with an exact sequence
$$0\seq{\cal O}\seq{\cal E}\seq{\cal I}_Y(4)\seq 0,$$
since by adjunction formula $\det N_{Y|{\bb P}^3}={\cal O}(4)|Y$.
If we choose concretely $Y$ as the image of $C:=\{y^2z-x^3+xz^2=0\}\subset{\bb P}^2$ via
{\small
\begin{eqnarray*}i:{\bb P}^2&\ratmap&{\bb P}^3\\ \
[x:y:z]&\mapsto & [x^2y+xyz+z^3:xy^2+yz^2+zx^2:x^2y+xyz+xz^2:xy^2+y^2z+z^3],
\end{eqnarray*}
}
we can compute with MACAULAY, that 
$$I_{Y}:=\bigoplus_{n\in{\bb N}}H^0({\cal I}_Y(n))$$
is generated by three quartics $q_1,q_2,q_3$ and four quintics. But also with MACAULAY we can verify that the projective schemes $Y$ and
$\{z\in{\bb P}^3|q_1(z)=q_2(z)=q_3(z)=0\}$ are identical. Hence ${\cal I}_Y(4)$ is generated by global sections and we conclude
that ${\cal E}$ is generated by global sections. Therefore $\can={\cal O}_{{\bb P}({\cal E})}(2)$ is globally generated and we can choose
$s\in H^0(\can)$ such that $X$ is smooth. 

By construction $c_1({\cal E})=4h$ and $c_2({\cal E})=8h^2$, hence $\gamma({\cal E})=-16$ and $|P|=128$.

\begin{remark}
\begin{enumerate}
\item Note that this case is extremal in some sense: Any elliptic curve in ${\bf P}^3$ of degree $d\ge 9$ cannot be cut out by
quartics. This can be seen like follows: If the contrary would be the case, the Serre construction would yield a globally generated
vector bundle ${\cal E}$ and a sequence
$$0\seq{\cal O}\seq{\cal E}\seq{\cal I}_Y(4)\seq 0.$$
Since $\can={\cal O}_{{\bb P}({\cal E})}(2)$ would then be globally generated as well, we conclude $(\can)^4\ge 0$. On the
other hand, $\gamma({\cal E})=16-4d\le -20$ and hence $(\can)^4=32\gamma+512\le -128$.
\item This example is extremal also in some other sense: A general member of $|\can|$ is an elliptic fibre space
over a quadric, which is the restriction of an elliptic fibre space ${\bf P}({\cal E})\seq{\bf P}^3$ (see \cite{diss}).
\end{enumerate}
\end{remark}
The cases where $\gamma$ is odd are more complicated to deal with. In these cases we cannot use elliptic curves. Instead of
genus $1$ we have to choose negative genera, hence $Y$ is not irreducible and the extendability
condition of the normal bundle is harder to check. So the Serre construction does not appear to be useful.


\begin{thebibliography}{Bot95}
\bibitem[AGV]{agv}Arnold, V.I., Gusein-Zade, S.M., Varchenko, A.N.: Singularities of Differential Maps Vol. II, Boston 1985-1988
\bibitem[Bo90]{bo90}Borcea, C.: Nodal quintic threefolds and nodal octic surfaces, Proc. of AMS 109, 627-635, 1990
\bibitem[Bot95]{bo} von Bothmer, H.-C.: Syzygien von glatten Raumkurven, Diplomarbeit, Bayreuth 1995
\bibitem[Ch92]{ch92}Chmutov, S.V.: Examples of projective surfaces with many singularities, J. Alg. Geom.
1, 191-196, 1992
\bibitem[Cl83]{cl83}Clemens, C. H.: Double Solids, Adv. in Math. 47, 107-230, 1983
\bibitem[CS99]{cs99}Cynk, S., Szemberg, T.: Double covers of ${\bb P}^3$ and Calabi-Yau varieties, math.AG/9902057, 1999
\bibitem[Cy99]{cy99}Cynk, S.: Double coverings of octic arrangements with isolated singularities, Adv. Theor. Math. Phys. 3, No.2, 217-225, 1999
\bibitem[En97]{en95}Endra{\ss}, S.: A projective surface of degree eight with 168 nodes, J. Algebr. Geom. 6, No.2, 325-334, 1997
\bibitem[Ga54]{ga54}Gallarati, D.: Sopra una superficie dell'ottavo ordine dotata di 157 nodi, Atti Accad. naz. Lincei, Rend., Cl. Sci. fis. Mat. natur., 
VIII. Ser. 16, 454-459, 1954 
\bibitem[Ga58]{ga58}Gallarati, D.: Una superficie dell'ottavo ordine con 160 nodi, Accad. Ligure Sci. Lett. 14, 1958
\bibitem[Ha]{ha}Hartshorne, R.: Algebraic Geometry, New York, 1977
\bibitem[Hu]{hu}Hulek, K.: Projective geometry of elliptic curves,
Ast\'erisque 137, Paris, 1986
\bibitem[K\"u01]{diss} K\"uhnel, M.: \"Uber gewisse Calabi-Yau-3-faltigkeiten mit Picardzahl $\rho(X)=2$,
doctoral thesis, Bayreuth 2001
\bibitem[Kr56]{kr56}Kreiss, H. O.: \"Uber syzygetische Fl\"achen, Annali die Matematica Pura et Appl. 41, 105-111, 1956
\bibitem[La]{labs}Labs, O.: http://www.algebraicsurface.net/octics/, 2003
\bibitem[Mi84]{mi84}Miyaoka, Y.: The maximal number of quotient singularities on surfaces with given numerical invariants, Mathematische Annalen 268, 159-171, 1984
\bibitem[OSS]{oss} Okonek, Ch., Schneider, M., Spindler H.: Vector bundles on Complex Projective Spaces, Boston (1988)
\bibitem[Sa01]{sa} Sarti, A:. Pencils of symmetric surfaces in ${\Bbb P}\sb 3$,  J. Algebra  246, no. 1, 429-452, 2001
\bibitem[Se53]{se53}Segre, B.: Sul massimo numero di nodi delle superficie algebriche, Atti Accad. Ligure Sci. Lett. 9, 15-22, 1953
\bibitem[We87]{we87}Werner, J.: Kleine Aufl\"osungen spezieller dreidimensionaler Variet\"aten, Bonner Mathematische Schriften 186, 1987
\end{thebibliography}
\end{document}